\documentclass[10pt, conference, compsocconf]{IEEEtran}
\usepackage[sort,compress]{cite}
\usepackage[margin=0pt,labelsep=space,scriptsize,labelfont=bf,textfont=it,up,belowskip=-8pt,aboveskip=5pt]{caption}
\usepackage{graphics,epsfig,graphicx,float,subfigure,color}
\usepackage{algorithmic}

\usepackage[section]{algorithm}

\usepackage{etex}

\usepackage{amsmath,amssymb,amsbsy,amsfonts,latexsym}%amsthm
\usepackage{multicol,multirow}
\usepackage{hhline}
\usepackage{paralist}
\usepackage{array}
\usepackage{marvosym}
\usepackage{url}
\usepackage{boxedminipage}
\usepackage{wrapfig}
 \usepackage[plainpages=false, colorlinks=true,
   citecolor=blue, filecolor=black, linkcolor=blue,
   urlcolor=blue]{hyperref}
\usepackage{multirow}
\usepackage{threeparttable}
\usepackage{pdfpages}
\usepackage{xspace}
\usepackage{colortbl}

\usepackage{tabularx,booktabs}

\usepackage[usenames,dvipsnames]{xcolor}

\usepackage{tikz} 
\usetikzlibrary{arrows} 
\usetikzlibrary{patterns}  
\usepackage{pgfplots} 
\pgfplotsset{compat=1.3}

\let\labelindent\relax
\usepackage{enumitem}

\usepackage{calligra}
\DeclareMathAlphabet{\mathcalligra}{T1}{calligra}{m}{n}

\newcommand{\figref}[1]{Figure~\ref{#1}}
\newcommand{\tabref}[1]{Table~\ref{#1}}
\newcommand{\secref}[1]{\S\ref{#1}}
\newcommand{\algref}[1]{Algorithm~\ref{#1}}

\renewcommand{\b}[1]{{#1}}

\newcommand{\MA}[1]{{\mathcal #1}}
\DeclareMathOperator{\bigO}{\mathcal{O}}

\definecolor{light-gray}{gray}{0.80}

\newcommand{\algcmt}[1]{\hfill {\footnotesize\ttfamily\textcolor{blue}{/* {#1} */}}}

\newcommand{\bu}{{\b u}} % pointential vector
\newcommand{\bw}{{\b w}} % weight vector
\newcommand{\bx}{{\b x}} % a point
\newcommand{\Ker}{\mathcal{K}}
\newcommand{\XX}{{\MA{X}}} % set of all points
\newcommand{\sk}[1]{{\widetilde{#1}}} % skeleton
\newcommand{\K}{{\MA{K}}} % The kernel matrix
\newcommand{\lc}{{\tt l}}  % left child of node $\alpha$
\newcommand{\rc}{{\tt r}}  % right child of node $\alpha$
\newcommand{\ns}{{s}} 
\newcommand{\IASKIT}{{\sc Inv-ASKIT}}
\newcommand{\ASKIT}{{\sc ASKIT}}
\newcommand{\GMRES}{{\sc GMRES}}

\expandafter\def\expandafter\normalsize\expandafter{
  \setlength\abovedisplayskip{6pt}
  \setlength\belowdisplayskip{6pt}
  \setlength\abovedisplayshortskip{6pt}
  \setlength\belowdisplayshortskip{6pt}
}

\renewcommand{\baselinestretch}{1.04}

\begin{document}

\title{\IASKIT{}: A Parallel Fast Direct Solver for Kernel Matrices}

\author{
\IEEEauthorblockN{Chenhan D. Yu\IEEEauthorrefmark{1}, William B. March\IEEEauthorrefmark{2},
Bo Xiao\IEEEauthorrefmark{3} and George Biros\IEEEauthorrefmark{4}}
\IEEEauthorblockA{
\IEEEauthorrefmark{1}Department of Computer Science\\
\IEEEauthorrefmark{1}\IEEEauthorrefmark{2}\IEEEauthorrefmark{3}\IEEEauthorrefmark{4}
Institute for Computational Engineering and Science\\
The University of Texas at Austin, Austin, Texas, USA \\
\IEEEauthorrefmark{1}chenhan@cs.utexas.edu,
\IEEEauthorrefmark{2}march@ices.utexas.edu,
\IEEEauthorrefmark{3}bo@ices.utexas.edu,
\IEEEauthorrefmark{2}gbiros@acm.org
}
}

\maketitle
We present a parallel algorithm for computing the approximate
factorization of an $N$-by-$N$ kernel matrix. Once this factorization
has been constructed (with $N \log^2 N$ work), we can solve linear
systems with this matrix with $N \log N$ work.  Kernel matrices
represent pairwise interactions of points in metric spaces.  They
appear in machine learning, approximation theory, and computational
physics. Kernel matrices are typically dense (matrix multiplication
scales quadratically with $N$) and ill-conditioned (solves can require
100s of Krylov iterations). Thus, fast algorithms for matrix
multiplication and factorization are critical for scalability.

Recently we introduced \ASKIT{}, a new method for approximating a
kernel matrix that resembles N-body methods. Here we
introduce \IASKIT{}, a factorization scheme based on \ASKIT{}. We
describe the new method, derive complexity estimates, and conduct an
empirical study of its accuracy and scalability.  We report results on
real-world datasets including ``COVTYPE'' ($0.5$M points in 54
dimensions), ``SUSY'' ($4.5$M points in 8 dimensions) and ``MNIST''
(2M points in 784 dimensions) using shared and distributed memory
parallelism. In our largest run we approximately factorize a dense
matrix of size 32M $\times$ 32M (generated from points in 64
dimensions) on 4,096 Sandy-Bridge cores. To our knowledge these
results improve the state of the art by several orders of magnitude.

\section{Introduction} \label{s:intro} Given a set of $N$ points $\bx_i \in \mathbb{R}^{d}$ and a kernel
function
$\Ker(\bx_i, \bx_j):\mathbb{R}^d\times\mathbb{R}^d\rightarrow \mathbb{R}$
, a kernel matrix is the $N\times N$ matrix%
\footnote{We restrict our
discussion to the case in which $\bx_i$ and $\bx_j$ come from the same
point set. The $i$-indexed points are called the target points, and
the $j$-indexed points are called the source points. A kernel matrix
can also involve target and source points that come from different
sets.}
 whose entries are given by $K_{ij} = \Ker(\bx_i,\bx_j)$ for
$i,j=1,\ldots,N$.  We are interested in the factorization of $K$.

Kernel matrices are typically dense, so matrix-vector multiplication
with $K$ (henceforth, \emph{\sc ``MatVec''}) requires $\bigO(N^2)$
work.  Solving a linear system with $K$ (which we will refer to
as \emph{\sc ``Invert'' K}) is not stable because $K$ can be nearly
singular. We typically invert $\lambda I +K$, where $\lambda$ is
the \emph{regularization parameter} and $I$ is the $N\times N$
identity matrix. The LU factorization of $\lambda I + K$ requires
$\bigO(N^3)$ work.  An iterative Krylov method requires $\bigO(N^2)$
work but the constant can be quite high depending on $\lambda$.

\textbf{Motivation and significance.}
Kernel matrices appear in N-body problems~\cite{greengard94}, machine
learning methods for classification, regression, and density
estimation \cite{gray-moore01, wasserman04,hofmann-scholkopf-smola08},
spatial statistics \cite{chen-anitescu14}, Gaussian processes
\cite{rasmussen-williams06}, and approximation
theory~\cite{schaback06}, among other applications. Typical operations
required with kernel matrices are matrix-vector multiplication,
solving linear systems, and eigenvalue computations. For example, in
machine learning, solving linear systems is used in ridge regression
and support vector machines. Eigenvalue computations are required in
spectral clustering. Such operations can be done either with iterative
solvers in which we only require matrix-vector multiplication.
However, the conditioning can be bad and convergence can be quite
slow.

For many (but not all) kernel matrices, dramatic acceleration of the
{\sc MatVec} and the {\sc Invert} operations can be achieved if we
exploit subtle properties of their structure.  For example, consider
the Gaussian kernel, a popular kernel in machine learning,
\begin{equation}\label{e:gaussian}
\Ker(\bx_i,\bx_j) = \exp\left(-\frac{1}{2} \frac{\|\bx_i-\bx_j\|_2^2}{h^2}\right),
\end{equation}
where $h$ is the {\em bandwidth}.  For small $h$, $K$ approaches the
 identity matrix whereas for large $h$, $K$ approaches the rank-one
 constant matrix. The first regime suggests sparse approximations
 while the second regime suggests global (numerical) low-rank
 approximations. For the majority of $h$ values, however, $K$ is
 neither sparse nor globally low-rank.  It turns out that
 the \emph{off-diagonal blocks} of $K$ may admit good low-rank
 approximations. A key challenge  is identifying and constructing
 these approximations. When the points are in low dimensions (say $d <
 10$), effective techniques for the {\sc MatVec} and {\sc Invert}
 operations exist.  But for high-dimensional problems, there is very
 little work.

{\bf Contributions:} Our goal is a method that exhibits algorithmic
and parallel scalability for inverting $\lambda I + K$ for large $N$
and $d$. We only require the ability to compute an entry of $K$ in
$\bigO(d)$ time.  Our scheme extends \ASKIT{}%
\footnote{\small
\emph{A}pproximate \emph{S}keletonization \emph{K}ernel \emph{I}ndependent
\emph{T}reecode. The \ASKIT{} library is available at~\url{http://padas.ices.utexas.edu/libaskit}.}%
, a method  we introduced in a series of
papers~\cite{march-xiao-biros15,march-xiao-yu-biros15,march-xiao-biros-fmm-e15}. \ASKIT{}
approximates $K$ in $\bigO( d N\log N)$ time. Here we introduce
\IASKIT{} for factorizing and inverting the \ASKIT{}-based
approximation of $K$. In particular:
\begin{itemize}[leftmargin=*]
\item We present sequential and parallel factorization algorithms that
have $\bigO(N\log^2 N)$ complexity and can be applied in $\bigO(N\log
N)$ time (described in \secref{s:seq}, \secref{s:par}).
\item We present a theoretical analysis of the complexity of the algorithm
      (in terms of time, communication and storage). The estimates for
      the parallel factorization and solution are given
      by \eqref{e:t-fact} and \eqref{e:t-solv} respectively
      in~\secref{s:par}.
\item An \ASKIT{} variant approximates $K$ as a triple sum:
a block-diagonal plus a low-rank plus a \emph{sparse} matrix. In this
case we use \IASKIT{} as a preconditioner for a GMRES iterative
solver.
\item We apply our solver to  a binary
      classification problem using kernel regression. We present
  the solver's accuracy, training and test errors, and compare it with
  Krylov iterative methods (described in \secref{s:results}).
\item We conduct strong and weak scaling analysis in which we report
  communication, computation and overall efficiency (described
  in \secref{s:results}).
\end{itemize}
Our method uses the hierarchical decomposition of \ASKIT{}. The
fundamental concept behind it is simple: First, approximate $K$ as the
sum of a block-diagonal matrix and a low-rank matrix. Second, compute
the inverse using the Sherman-Morrison-Woodbury (SMW) formula. Third,
apply this decomposition recursively to each block of the
block-diagonal matrix.  \emph{To our knowledge, \IASKIT{} is the only
algorithm that can be used with high-dimensional data for a large
variety of kernels and the only algorithm that supports shared and
distributed memory parallelism. \IASKIT{} does not assume symmetry of
the kernel, global low-rank, sparsity, or any other property other
than that, up to a small sparse matrix-correction, the off-diagonal
blocks admit a low-rank approximation.}

{\bf Related work:} We do not review the large literature of
accelerating the matrix-vector operation with a kernel matrix. We
focus only on fast factorization schemes.  Many times kernel matrices
admit a good \emph{global low-rank} approximation.  For such matrices
the \emph{Nystrom method} and its variants can be very
effective~\cite{williams-seeger01,gittens-mahoney13,halko-martinsson-tropp11}
for both approximating $K$ and its regularized inverse.  However, in
many applications $K$ does not have a global low-rank structure and
Nystrom methods fail. In~\cite{march-xiao-biros-e15}, we discussed
this in detail and presented classification results on real datasets
for which the Nystrom method fails.

The idea of exploiting the block-diagonal-plus-low-rank structure has
been discussed extensively in the past fifteen years for kernel
matrices in low dimensions and finite-difference/finite element
matrices. Indeed for point datasets in low dimensions, there exist
methods that are much faster and more scalable
than \ASKIT{}. Approximating $(\lambda I + K)^{-1}$ is more
challenging, but this is also a well-studied problem.  Efficient
schemes for approximating the inverse of kernel (but also other types
of) matrices are also known as \emph{fast direct solvers}, and they
have roots in nested dissection (or domain decomposition),
multifrontal methods and truncated Newton-Schulz iterations.  Examples
include~\cite{ambikasaran-darve13,ambikasaran-13,greengard-neil-e14,ho-greengard12,greengard-gueyffier-martinsson-rokhlin09,bebendorf08,wang-li-dehoop13,le2007domain,hackbusch2008approximate}
and this list is non-exhaustive.  All these works, however, have been
applied successfully only in low-dimensional problems ($d<4$) and do
not address parallelism and performance.

To our knowledge, \cite{kriemann2005parallel} is the earliest work
that parallelizes hierarchical matrices with bulk synchronous
parallelization on shared memory systems, but it is not efficient due to
the $\MA{O}(N)$ critical path. In~\cite{66699}, $\MA{H}$-LU improves
its efficiency by employing a DAG-based task scheduling.  Distributed
implementations include 
\cite{wang-li-dehoop13} and \cite{izadi2012parallel}, but, to our
understanding, they were designed for sparse matrices related to 2D
and 3D stencil discretizations, not kernel matrices.  Although the
basic machinery is similar, but the matrices are quite different. We
discuss this more in \secref{s:hier}.  For Gaussian kernels, the only
work for fast inversion
is~\cite{greengard-neil-e14,ambikasaran-darve13} in which the authors
consider Gaussian and other radial-basis kernels in low dimensions
($d<4$). Parallelism and scalability were not considered in that
work. Their largest test had 1M points in 3D. We consider datasets
that have 30$\times$ more points and over 200$\times$ higher
dimensions.

{\bf Limitations:} This is a novel algorithm and this paper represents
our first attempt in an efficient implementation. For this reason it
has several limitations. First, we do not provide a theoretical
analysis of its numerical stability. (But it is relatively easy to
prove it for symmetric and positive definite matrices.) Second, we do
not exploit certain symmetries in \ASKIT{} that can further reduce the
factorization time. For example, $K$ can be symmetric. Third, we do
not performance-tune various algorithm parameters. Fourth, we do not
discuss optimizations that can be used to amortize the factorization
cost over different values of $\lambda$. Fifth, for certain \ASKIT{}
variants our scheme can be only used as a preconditioner and the
errors can be significant. Sixth, there are additional low-rank
structures that we do not exploit.  Finally, there exist fundamental
limitations inherited from \ASKIT{}. Not all kernel matrices have a
block-diagonal-plus-low-rank structure. Examples are point datasets
that have high \emph{ambient dimensionality} and kernels that are
highly oscillatory (e.g., the high-frequency Helmholtz kernel). Also,
the sparse-matrix correction for the off-diagonals can have large
norm, in which case using this version of \IASKIT{} as a preconditioner
will not work.

\section{Methods} \label{s:methods} \begin{table}
\centering
{\footnotesize
  \begin{tabular}{lll}
  \toprule
  Notation        &
                  & Description \\
  \cline{1-1} \cline{3-3} 
  $N$             &
                  & problem size \\
  $d$             &
                  & dimension of input points\\
  $m$             & 
                  & leaf node size \\
  $s$             & 
                  & number of skeleton points \\
%  $\kappa$        &
%                  & neighbor size \\
%  $\tau$          &
%                  & tolerance in adaptive skeletonization \\                  
%  $l$             &
%                  & tree level \\
  $p$              &
                   & number of MPI processes \\
%  $n$             &
%                  & number of points per MPI process \\
%  $r$             &
%                  & number of right hand sides \\                                 
  $\alpha$, $\beta$ &
                  & tree nodes, overloaded with the indices they own   \\               
  $\lc$, $\rc$ &
                  & left and right children of the treenode \\                  
  \toprule
  $\XX$           & 
                  & $\mathbb{R}^{d \times N}$, point dataset \\
  $\XX(\alpha)$   & 
                  & points owned by node $\alpha$, i.e.  $\XX(\alpha)=\{ \bx_i \lvert \forall i \in \alpha\}$ \\
  $\sk{K}$       &
                  & approximate kernel matrix  \\
  \toprule
  \end{tabular}
}
  \caption{Main notation used. Also occasionally we will use MATLAB-style matrix 
    representations.}
\label{tab:notation}
\end{table}
\textbf{Notation:}
In \tabref{tab:notation}, we summarize the main notation used in the paper. 
%We use $\alpha$ and $\beta$ to denote two siblings in the binary tree. 
We begin with discussing the basic concepts in factorization of
hierarchical matrices~\secref{s:hier}. Then we present a brief summary
of the \ASKIT{} algorithm~\secref{s:askit}.  We conclude with the
sequential~\secref{s:seq} and parallel~\secref{s:par} factorization
algorithms.

  \subsection{Hierarchical Matrices and Treecodes} \label{s:hier} Informally, $K\in \mathbb{R}^{N\times N}$ is hierarchical if it can be
approximated by $D + UV$, where $D$ is a \textbf{\emph{
block-diagonal matrix}}, $UV$ is \textbf{\emph{low-rank}},
and \textbf{\emph{each block of $D$ is also a hierarchical
matrix}}. To define the $D$, $U$, $V$ matrices first we consider the
following block partitioning to $N/2$-by-$N/2$ blocks
\begin{equation}
  \label{e:partitioning}
  K =
  \begin{bmatrix}
  K_{11} & K_{12}\\
  K_{21} & K_{22} \\
  \end{bmatrix} \\
   = 
\begin{bmatrix}
K_{11} & 0 \\ 0 & K_{22} \\ \end{bmatrix} + \begin{bmatrix} 0 &
K_{12} \\ K_{21} & 0 \\ \end{bmatrix}.  \\ \end{equation} $D$ is the
first matrix in the sum and the product of $UV$ is equal to the second
matrix; $U$ and $V$ will be defined shortly.  So in summary: (1)
Blocks $K_{12}$ and $K_{21}$ admit low-rank \emph{approximations},
which implies that $U$ and $V$ are low-rank.  (2) The same
decomposition applies recursively for the diagonal blocks $K_{11}$ and
$K_{22}$.  This approximation derived by recursive partitioning is
directly associated with a binary partitioning of the rows and columns
of $K$.  Note that using lexicographic ordering of the rows like
in \eqref{e:partitioning} will not, in general, result in low-rank
off-diagonal blocks. An appropriate reordering, in our case using the
tree structure, is critical for exposing such low-rank structure.
Given the finest granularity $m$, we partition the points using
geometric information and a binary tree so that we have less than $m$
points per leaf.  (The number of points $m$ is manually
selected. Typically is about 1048 or 2048 and it effects both
computation time and accuracy. For a discussion,
see~\cite{march-xiao-biros15}).  We denote the levels of the tree by
$l$, with $l=0$ at the root and $l=\log_2(N/m)$ at the leaves. Let
$\alpha$ be a node in the tree. We ``overload'' the notation by using
$\alpha$ to also indicate the indices of the points assigned to
$\alpha$ so that all $\XX(\alpha)
\in \alpha$. Then the hierarchical partition and approximation of $K$
assumes that for every node $\alpha$, we can partition the
self-interactions between its points as follows
\begin{equation}
  \label{e:partitioning-node}
  K_{\alpha\alpha} = 
  \begin{bmatrix}
  K_{\lc\lc} & K_{\lc\rc} \\
  K_{\rc\lc} & K_{\rc\rc} \\
  \end{bmatrix} = 
  \begin{bmatrix}
  K_{\lc\lc} & \\
  & K_{\rc\rc} \\
  \end{bmatrix} +
  \begin{bmatrix}
  & K_{\lc\rc} \\
  K_{\rc\lc} & \\
  \end{bmatrix}, 
\end{equation}
where $\lc$ denotes the left child of a node, and $\rc$ the right
child of a node. Equivalently to \eqref{e:partitioning}, we can write
$\sk{K}_{\alpha \alpha} = D_{\alpha} + U_{\alpha}V_{\alpha}$.  The
matrices $U_\alpha$ and $V_\alpha$ are selected to approximate
off-diagonal blocks in \eqref{e:partitioning-node} such that
\begin{equation}\label{e:uv}
  \begin{bmatrix}
  0 & K_{\lc\rc} \\
  K_{\rc\lc} & 0 \\
  \end{bmatrix} \approx
  U_\alpha V_\alpha.
\end{equation}
The particular form of these matrices for the \ASKIT{} method is given
in \eqref{e:uv}. If $K$ is invertible, we can write immediately the
\textbf{\emph{approximate inverse}} of $K$ by using the SMW formula:%
%
%\footnote{Here we assume $K$ is invertible. The extension to $\lambda I + K$ is trivial.}
%
\begin{equation}\label{e:smw}
\begin{split}
\sk{K}   &= D+UV = D(I+D^{-1}U V) \\
&= D(I+WV), \mbox{~where~} W=D^{-1}U.\\
\sk{K}^{-1} &= (I+WV)^{-1} D^{-1}\\
      &= (I-W(I+VW)^{-1}V)D^{-1} \\
      &= (I-WZV)D^{-1}, \mbox{~where~} Z=(I+VW)^{-1}.
\end{split}
\end{equation}
The extension to $\lambda I + K$ is trivial.
In \IASKIT{}, factorizing $K$ amounts to visiting every tree node
$\alpha$ to compute $W$ and $Z$. We require traversing the descendants
of $\alpha$ to compute $D^{-1}U$. If $\ns$ is the numerical rank of
$K_{\lc\rc}$ and $K_{\rc\lc}$, then computing $W$ requires $2\ns$
linear solves with $D_\alpha$. Computing $Z$ requires factorizing
(using LAPACK) the $2\ns\times 2\ns$ matrix $I+VW$.  The recursion
stops at the leaves of the tree in which we have a single matrix $
D_\alpha=K_{\alpha\alpha}$, which we then again factorize using
LAPACK.  We give pseudocode for the precise algorithms in~\secref{s:seq}.

As mentioned in \secref{s:intro}, there is extensive work on fast
solvers for hierarchical matrices. Such solvers rely on the SMW
formula. But the details can differ significantly, for example the
permutation or reordering of the matrix, construction of $U$ and $V$,
and the exploitation of various symmetries.%
\footnote{More precisely our matrix can be described as a HODLR
  (Hierarchical Off-Diagonal Low Rank)
  matrix~\cite{ambikasaran-13}. Other types of matrices include hierarchically
  semi-separable matrices (HSS), and
  $\MA{H}$-matrices \cite{hackbusch1999sparse}.} 
\ASKIT{} can also produce a slightly different
decomposition (known as a fast multipole matrix or an $\MA{H}$-matrix), 
which decomposes $K$ into a sum of a hierarchical matrix
$D+UV$ (which has the form we discussed above) and
a \textbf{\emph{sparse matrix}} $S$.  Next, we give some details
for \ASKIT{}'s construction of the hierarchical approximation of $K$
and briefly explain the $\sk{K}= D+UV+S$ decomposition.

  \subsection{ASKIT} \label{s:askit} %%%%%%%%%%%%%%%%%%%%%%%%%%%%%%%%%%%%%%%%%%%%%%%%%%%%%
Details on the algorithm can be found
in~\cite{march-xiao-biros15}. Here we briefly summarize it. The input
is a set of points $\XX$ and a list of nearest neighbor points in
$\XX$ for each point $\bx_i \in \XX$. We build a balanced binary
tree-based partitioning of the points. We use a parallel ball
tree~\cite{omohundro1989five} data structure. For $d<4$ quadtrees and
octrees can be effectively used, but in high dimensions such data
structures cannot be used.

After building the tree, we construct the hierarchical representation.
In \ASKIT{} we use the notion of a point being \emph{separated} from a
tree node $\alpha$. A basic separation criterion is that
$i \notin \alpha$.  If $K_{i\alpha}$ is the interaction
between $i$ and $\alpha$, we construct $K_{i\sk{\alpha}}$ and
$P_{\sk{\alpha} \alpha}$ such that
\begin{equation}\label{e:askit-point}
K_{i\alpha} \approx K_{i\sk{\alpha}} P_{\sk{\alpha} \alpha}.
\end{equation}
The $\sk{\alpha}\subset \alpha$ points are referred to as the
``\textbf{\emph{skeletons}}''.  The skeletons are the pivot columns of
the QR factorization of a subsample of the block matrix
$K_{[1:N]\alpha}$. (The nearest-neighbor information is used to
effectively sample this block matrix.)
$K_{i\sk{\alpha}}=\Ker(\bx_i,\XX(\alpha))$ does not involve any
approximation, but instead of computing interactions with all points
in $\alpha$ we just compute interactions with
$\sk{\alpha}$. $P_{\sk{\alpha} \alpha}$ is computed using the
interpolative decomposition algorithm.  A key property of \ASKIT{} is
that the skeletons of a node are selected from the skeletons of its
children. That is, for a node $\alpha$ the skeletons $\sk{\alpha}$ are
a subset of $\sk{\lc} \cup \sk{\rc}$ or
$[\sk{\lc} \sk{\rc}]$. Furthermore
$P_{\sk{\alpha} \alpha} \bw(\alpha)$ (when we compute $K\bw$) can be
computed recursively by only storing $P_{\sk{\alpha}
[\sk{\lc}\sk{\rc}]}$. If we assume the skeleton size $\ns$ is fixed
for all nodes, then $P_{\sk{\alpha} [\sk{\lc}\sk{\rc}]}$ is an
$\ns \times 2\ns$ matrix.

The details of constructing skeletons and the $P$ matrices can be
found in~\cite{march-xiao-biros15}.  In the pseudocode below we outline
the two steps required to perform a fast {\tt MatVec} $\bu=K\bw$.  The
``\emph{Skeletonize}'' step builds the projection matrices $P$ for each
node and computes $\bw_{\sk{\alpha}} =
P_{\sk{\alpha} \alpha}\bw_\alpha$, auxiliary variables that are then
used in the ``\emph{Evaluate}'' step in which for every target point
$i$, the tree is traversed top-down to compute $\bu_i$ the entry of
the output vector for point $i$.
%\footnote{The algorithm below shows
%evaluation for a single target point $i$. A simple ``for loop''
%computes the solution for all target points.}

\
\hrule

\begin{minipage}[l]{0.42\linewidth}
\small
\begin{algorithmic}[1]
    \STATE \emph{//Skeletonize:}
    \STATE {\tt Skeletonize}$(\alpha)$
  \IF {$\alpha$ is leaf}
    \STATE Find $\sk{\alpha}$
    \STATE Construct $P_{\sk{\alpha} \alpha}$
    \STATE $\bw_\sk{\alpha} = P_{\sk{\alpha} \alpha} \bw_\alpha$ 
  \ELSE
    \STATE {\tt Skeletonize}$(\lc)$
    \STATE {\tt Skeletonize}$(\rc)$
    \STATE Construct $\sk{\alpha}$, $P_{\sk{\alpha}[\sk{\lc} \sk{\rc}]}$
    \STATE $\bw_\sk{\alpha} = P_{\sk{\alpha} [\sk{\lc} \sk{\rc}]}  \bw_{[\sk{\lc} \sk{\rc}]}$ 
  \ENDIF
\end{algorithmic}
\end{minipage}
\hfill
\begin{minipage}[r]{0.57\linewidth}
\small
\begin{algorithmic}[1]
  \STATE \emph{//Evaluate:}
  \STATE {\tt MatVec} $(i,\alpha)$
  \IF {{\tt Separated}$(i,\alpha)$}
     \STATE $\bu_i = K_{i\sk{\alpha}} \bw_{\sk{\alpha}}$
  \ELSE
     \IF {$\alpha$  is leaf}
        \STATE $\bu(i)+=K_{i \alpha} \bw_\alpha$
     \ELSE 
       \STATE {\tt MatVec}$(i,\lc)$
       \STATE {\tt MatVec}$(i,\rc)$
     \ENDIF
  \ENDIF
\end{algorithmic}
\end{minipage}
\vspace*{7pt}

\hrule
\

In the \emph{Skeletonize} algorithm, lines 4, 5 and 10 do not depend
on $\bw$. Therefore, if we want to apply $K$ to a different vector
$\bw$, only lines 6 and 11 are needed. In the \emph{Evaluate}
phase, \ASKIT{} provides two different algorithms for the {\tt
Separated$(i,\alpha)$} function.  In the first  algorithm,
$i$ is separated from $\alpha$ if $i \notin \alpha$. This definition
creates a hierarchical approximation of $K$ identical to
 \eqref{e:partitioning-node}.  Then, if $\lc$ and $\rc$ are two sibling
nodes with parent $\alpha$, the matrices $K_{\lc\rc}$ and
$K_{\rc \lc}$ are approximated by
\begin{equation}\label{e:offdiag-askit}
K_{\lc \rc} \approx K_{\lc \sk{\rc}} P_{\sk{\rc} \rc} \mbox{~and~}
K_{\rc \lc} \approx K_{\rc \sk{\lc}} P_{\sk{\lc} \lc}.
\end{equation} 
We now can define for a node $\alpha$
\begin{equation}\label{e:uv-askit}
U_\alpha  = 
\begin{bmatrix}
\K_{\lc\sk{\rc}} & 0 \\
0 & K_{\rc\sk{\lc}}
\end{bmatrix}
\mbox{~and~}
V_\alpha =
\begin{bmatrix}
0 & P_{\sk{\rc}\rc}\\
P_{\sk{\lc}\lc} & 0
\end{bmatrix}.
\end{equation}
This completes the definition of the hierarchical matrix for \IASKIT{}
(compare \eqref{e:uv-askit} with \eqref{e:uv}).  To reduce memory costs
however, \ASKIT{} never constructs $P_{\sk{\alpha} \alpha}$
explicitly. Instead, looking at the \emph{Skeletonize}
function, \ASKIT{} only stores $P_{\sk{\alpha} [\sk{\lc}\sk{\rc}]}$
for each node $\alpha$. In the factorization we will need the whole
$P_{\sk{\alpha}\alpha}$ matrix. We explain how to compute it
in~\secref{s:seq}.

Other then $i \notin \alpha$, \ASKIT{} can also define $i$ to be
separated from $\alpha$ 
%if $\alpha$ does not contain any \emph{nearest-neighbors} of $i$. 
if no \emph{nearest-neighbor} of $i$ is in $\alpha$.
This definition introduces a
more complicated structure in the matrix. We do
construct \eqref{e:offdiag-askit} but for points $i$ that have 
\emph{nearest-neighbors} in $\alpha$, line 4 in the \emph{Evaluate} is replaced by
$\bu_i = K_{i\alpha}\bw_\alpha$, i.e., an exact evaluation. For
certain kernels this significantly increases the accuracy of the
overall scheme. A convenient way to think about this is that the
second algorithm builds an approximation that is equivalent to
\begin{equation}\label{e:askit}
K \approx (D + UV) + S = \sk{K}+S,
\end{equation}
where $D+UV$ is the hierarchical matrix and $S$ is a sparse matrix
that is zero if we use the $i\notin \alpha$ separation criterion.  If
the norm of $S$ is small relative to $\sk{K}$ we expect that $\sk{K}$
will be a good preconditioner for $\sk{K}+S$. (Again there are kernels
for which this is not the case.)

In this section we summarized the basic ideas in \ASKIT{}. Its
structure is very similar to what is known as a \emph{treecode} in
computational physics. We omitted several important details of
the \ASKIT{} algorithm. For example, we assumed that the number of
skeleton points $\ns$ is fixed. In reality, this number differs for
each tree node and it is chosen adaptively so that the approximation
error is below a tolerance $\tau$. Finally, the parallelization and
scalability of \ASKIT{} was described
in~\cite{march-xiao-yu-biros15}. Another version of \ASKIT{} that
resembles the relationship between fast-multipole methods and
treecodes is described in~\cite{march-xiao-biros-fmm-e15}. In these
papers, we also discuss the error analysis and work/storage complexity
of \ASKIT{}. Next, we will assume that we are given $D+UV$
from \ASKIT{} and we proceed with discussing the \IASKIT{}
factorization in detail using the ideas we outlined
in \secref{s:hier}.

  \subsection{Fast direct solver} \label{s:seq} Given the binary tree, the regularization parameter $\lambda$, and
$\sk{K}$ in terms of $\sk{\alpha}$ and $P_{\sk{\alpha}\alpha}$ 
for every tree node $\alpha$, 
we construct a factorization of $\lambda I +\sk{K}$ that
allows the fast solution of linear systems.  We describe the algorithm
for the case of $\lambda=0$, for notational simplicity, but the
algorithms are stated for any $\lambda$. So given
$\bu\in\mathbb{R}^{d}$, we wish to compute $\bw = \sk{K}^{-1}\bu$.

%\mbox{, where } \sk{K}\approx\lambda I + K.
%\end{equation}
In the factorization phase we compute $U$, $W$, $V$ and $Z$.  We first
show how \eqref{e:smw} is computed in a blocked $2\times2$
matrix. $D(I+WV)$ is computed as
\begin{equation}
\label{e:hodlr}  
  \begin{bmatrix}
  \sk{K}_{\lc\lc} &  \\
  & \sk{K}_{\rc\rc} \\
  \end{bmatrix}
  \left(
  I
  +
  \begin{bmatrix}
  W_{\lc\sk{\rc}} &  \\
  & W_{\rc\sk{\lc}} \\
  \end{bmatrix}
  \begin{bmatrix}
  & P_{\sk{\rc}\rc} \\
  P_{\sk{\lc}\lc} & \\
  \end{bmatrix}
  \right).
\end{equation}
where $W_{\lc\sk{\rc}}=\sk{K}_{\lc\lc}^{-1}K_{\lc\sk{\rc}}$
and $W_{\rc\sk{\lc}}=\sk{K}_{\rc\rc}^{-1}K_{\rc\sk{\lc}}$.
Then the partitioned form of $I-WZV$ is 
\begin{equation}
\label{e:blocksmw}
  I - 
  \begin{bmatrix}
  W_{\lc\sk{\rc}} &  \\
  & W_{\rc\sk{\lc}} \\
  \end{bmatrix}
  \begin{bmatrix}
  I & P_{\sk{\rc}\rc}W_{\rc\sk{\lc}} \\ 
  P_{\sk{\lc}\lc}W_{\lc\sk{\rc}} & I \\
  \end{bmatrix}^{-1}
  \begin{bmatrix}
  & P_{\sk{\rc}\rc} \\
  P_{\sk{\lc}\lc} & \\
  \end{bmatrix}.
\end{equation}
%where $Z$ composes $P_{\sk{\lc}\rc}W_{\rc\sk{\lc}}$ 
%and $P_{\sk{\rc}\lc}W_{\lc\sk{\rc}}$.
If $\alpha$ is an non-leaf node, $P_{\sk{\alpha}\alpha}$ is evaluated  recursively by
\begin{equation}
  \label{e:innerid}
  P_{\sk{\alpha}\alpha} = 
  P_{\alpha[\sk{\lc}\sk{\rc}]}
  \begin{bmatrix}
  P_{\sk{\lc}\lc} &  \\
  & P_{\sk{\rc}\rc} \\
  \end{bmatrix}
  \mbox{ using \texttt{GEMM}.}
\end{equation} 
At the leaf nodes $P_{\sk{\alpha}\alpha}$ is computed by \ASKIT{}
(using the LAPACK routines \texttt{GEQP3}
and \texttt{GELS}~\cite{march-xiao-biros15}).

\begin{algorithm}[!htp]
\begin{algorithmic}
  \IF { $\alpha$ is leaf}
    \STATE LU factorization $\lambda I + K_{\alpha\alpha}$. \algcmt{$\MA{O}(m^3)$}
    %\STATE Compute $\SKa$ and \eqref{e:leafid}. \algcmt{$\MA{O}(ms^2)$}
  \ELSE
    \STATE \texttt{Factorize}($\lc$) and \texttt{Factorize}($\rc$). %\algcmt{$2T^f(\frac{N}{2})$}
    %\STATE Split out $E$ using $\kappa$-nearest neighbors.
    \STATE Compute $K_{\lc\sk{\rc}}$, $K_{\rc\sk{\lc}}$, $P_{\sk{\lc}\lc}$, $P_{\sk{\rc}\rc}$.  \algcmt{$\MA{O}(s^3+Ns^2)$}
    \STATE $W_{\lc\sk{\rc}}=\texttt{Solve}(\lc, K_{\lc\sk{\rc}})$ and $W_{\rc\sk{\lc}}=\texttt{Solve}(\rc, K_{\rc\sk{\lc}})$. %\algcmt{$2T^s(\frac{N}{2},s)$}
    %\STATE LU factor $Z$ in \eqref{e:blocksmw}. \algcmt{$O(s^3)$}
    \STATE LU factorize the reduced system in \eqref{e:blocksmw}. \algcmt{$O(s^3)$}
  \ENDIF
\end{algorithmic}
\caption{{} \texttt{Factorize}($\alpha$)}
\label{a:factor}
\end{algorithm}

\textbf{Factorization Phase.}
\algref{a:factor} performs a postorder traversal of the tree. 
If $\alpha$ is a leaf node, we factorize $\lambda I +
K_{\alpha\alpha}$ with partial-pivoting LU factorization (LAPACK
routine \texttt{GETRF}).  Otherwise, we first compute
$K_{\lc\sk{\rc}}$, $K_{\rc\sk{\lc}}$, $P_{\sk{\lc}\lc}$,
$P_{\sk{\rc}\rc}$ using \eqref{e:innerid}.  Then we solve
$W_{\lc\sk{\rc}}=\sk{K}_{\lc\lc}^{-1}K_{\lc\sk{\rc}}$ and
$W_{\rc\sk{\lc}}=\sk{K}_{\rc\rc}^{-1}K_{\rc\sk{\lc}}$ by calling
\texttt{Solve}() (\algref{a:solve}). (Note that when we visit a node at some level, the
factorization information required for the solve is already in place.)
We compute and factorize $Z$ using the blocked representation
in \eqref{e:blocksmw}.  \texttt{Solve}() is described
in \algref{a:solve}, which applies $\sk{K}_{\alpha\alpha}^{-1}$ to a vector $\bu$.
\begin{algorithm}[!htp]
\begin{algorithmic}
  \IF {$\alpha$ is leaf}
    \STATE LU solver $\bw = (\lambda I + K_{\alpha\alpha})^{-1}\bu$. \algcmt{$\MA{O}(m^2)$}
  \ELSE
    \STATE $v_{\lc}=\texttt{Solve}(\lc, \bu_{\lc})$ and $v_{\rc}=\texttt{Solve}(\rc, \bu_{\rc})$. %\algcmt{$2T_s(\frac{N}{2},r)$}
    \STATE $v = [v_\lc\ v_\rc]$
    \STATE Compute $\bw=v-WZVv$ using \eqref{e:blocksmw}. \algcmt{$\MA{O}(sN)$}
  \ENDIF
\end{algorithmic}
\caption{{} $\bw=\texttt{Solve}(\alpha,\bu)$} 
\label{a:solve}
\end{algorithm} 
 
\textbf{Solving Phase.}
If $\alpha$ is a leaf node, $\bw = K_{\alpha\alpha}^{-1}\bu$ directly
invokes an LU solver (LAPACK routine \texttt{GETRS}).  Otherwise, we first solve
$v_{\lc}=\sk{K}_{\lc\lc}^{-1}\bu_{\lc}$ and
$v_{\rc}=\sk{K}_{\rc\rc}^{-1}\bu_{\rc}$ recursively.  Then we
compute \eqref{e:blocksmw} with two matrix-matrix multiplications
(BLAS routine \texttt{GEMM}) and an LU solver (\texttt{GETRS}).  Note that both
algorithms are implemented using bottom-up tree traversal and ``for-
loops''---we do not actually employ recursion.

\textbf{Refactoring for different $\lambda$.}
In the case that we want to refactor for a different value of
$\lambda$, (e.g., in cross-validation studies), $\sk{K}$ needs to be
refactorized, but the skeletons and $V$ we have computed remain the
same. Thus, when we refactor $\sk{K}$ with different $\lambda$, only
$W$ and $Z$ need to be recomputed. (The factorization at the leaf
level can be also amortized if SVD is used instead of LU, 
but we do not discuss this here.)

%Despite that \algref{a:factor} is a recursive function, a for-loop bottom up
%level by level traversal is implemented to avoid recursion.
%Each node computes its skeletons and evaluate $U$, $V$ required by 
%it's parent
%using \eqref{e:leafid} or \eqref{e:innerid}.
%Leaf nodes factorize diagonal blocks with LU factorization, and non-leaf
%nodes factorize the reduce system from \eqref{e:smw} after solving
%$\sk{U}_1$ and $\sk{U}_2$.
%Overall, the complexity of factorization $T_f(N)$ is $\MA{O}(s^2N\log^{2}{N})$.
%In the case that $s$ is a constant related to $N$, the complexity can be
%simplified as $\MA{O}(N\log^{2}{N})$.

%We start from constructing the robust prepresentation. 
%Given a complete binary tree that nestedly partitions $\XX$ into leaf nodes,
%\algref{a:factor} constructs $K=H+E$ from bottom to the root of the tree.
%Each leaf node $\alpha$ own a block diagonal $K(\alpha,\alpha)$, and a 
%partial-poviting LU (\texttt{GETRF}) is computed and stored.
%Skeletons ($\SKa$) is computed by \algref{}, and the full projection matrix 
%($P_\alpha$) is parent's $V$. Parent's $U$ is computed with sibling's 
%skeletons.
%For an non-leaf node, skeletons are computed after the recusion, but the
%projection matrix is $s \times 2s$, which can't be directly used as $V$.
%To

%Overall, the complexity of \algref{a:solve} $T_s(N,nrhs)$ is $\MA{O}(nrhs
%    \times s N\log{N})$.

%This variant doesn't change the worst cast complexity, and we still
%use $s$ in the following complexity analysis.

\textbf{Complexity of the sequential algorithm.}
We present the complexity analysis of \algref{a:factor}
and \algref{a:solve}. The main parameters are $\ns$ and $m$. They
depend on the kernel, the dimensionality, and the geometric
distribution of points. They control the accuracy of the scheme, but
are not explicitly related to problem size $N$. In the actual
implementation $\ns$ is chosen adaptively for each node.  If we use its
maximum value (across nodes), then the results below are an upper bound to
the work complexity. Note that in some applications, the kernel does
depend on the problem size, for example the bandwidth value for the
Gaussian kernel in machine learning problems. In that case $m$ and
$\ns$ will need to change accordingly. But here we assume a fixed
kernel for all values of $N$, since the way the kernel changes with
increasing $N$ depends on the application and the geometry of the
dataset.

$T^f(N)$ denotes the complexity of the factorization \algref{a:factor}
with $N$ points, and $T^s(N)$ denotes the complexity of 
\algref{a:solve} with $N$ points. We can verify that
\begin{equation}
  T^s(N) = 2T^s\left(\frac{N}{2}\right)+\MA{O}(Ns) = \MA{O}(sN\log{N}).
\end{equation}
Using the solving phase complexity, we derive $T^f(N)$ as
\begin{equation}
  T^f(N) = 2T^f\left(\frac{N}{2}\right)+s^2T^s(N) = \MA{O}(s^2N\log^{2}{N}).
\end{equation}
To summarize, \algref{a:factor} takes $\MA{O}(N\log^{2}{N})$ work
and \algref{a:solve} takes $\MA{O}(N\log{N})$ work.
 
\textbf{Storage.}
Beyond the storage required for \ASKIT{}, we need to store $U$
(overwritten by $W$), $V$ and the LU factorization in each level.  $U$
and $V$ take $\MA{O}(2sN)$ in each level.  Overall, the space
requirement is $\MA{O}(dN+2sN\log{\frac{N}{m}})$.  To save about half
of the memory, $V$ can be constructed recursively when it's needed, and
we discard it when we reach the parent level. This will result in
slightly larger constants in the solving phase.

\textbf{Shared memory parallelism:}
There are several ways to exploit shared memory parallelism. We have
implemented and tested two of those.  The first way is level-by-level
parallel tree traversal for both construction and evaluation, in which
we use sequential versions of BLAS/LAPACK routines, concurrent
traversal of nodes at the same level, and sequential traversal across
levels. The second way is to use sequential traversals and use the
shared-memory parallelism for the BLAS/LAPACK routines. It turns out,
that it is more efficient to use a hybrid version of these two
approaches.  From level-$4$ to the leaf level, we disable parallel
BLAS and LAPACK because we have experimentally observed that parallel
tree traversal provides higher efficiency due to better memory
locality.  From level-$0$ to level-$3$ in the local tree, we switch
from parallel tree traversal to parallel BLAS operations. In parallel
tree traversal some cores may start to idle from tree level three,
which only has 8 tree nodes.  So in the top four levels of the tree
we use parallel \texttt{GETRF}, \texttt{GETRS} and \texttt{GEMM}.

 \subsection{Parallel algorithm} \label{s:par} 
\begin{figure}[!t]
  \centering \includegraphics[scale=.4]{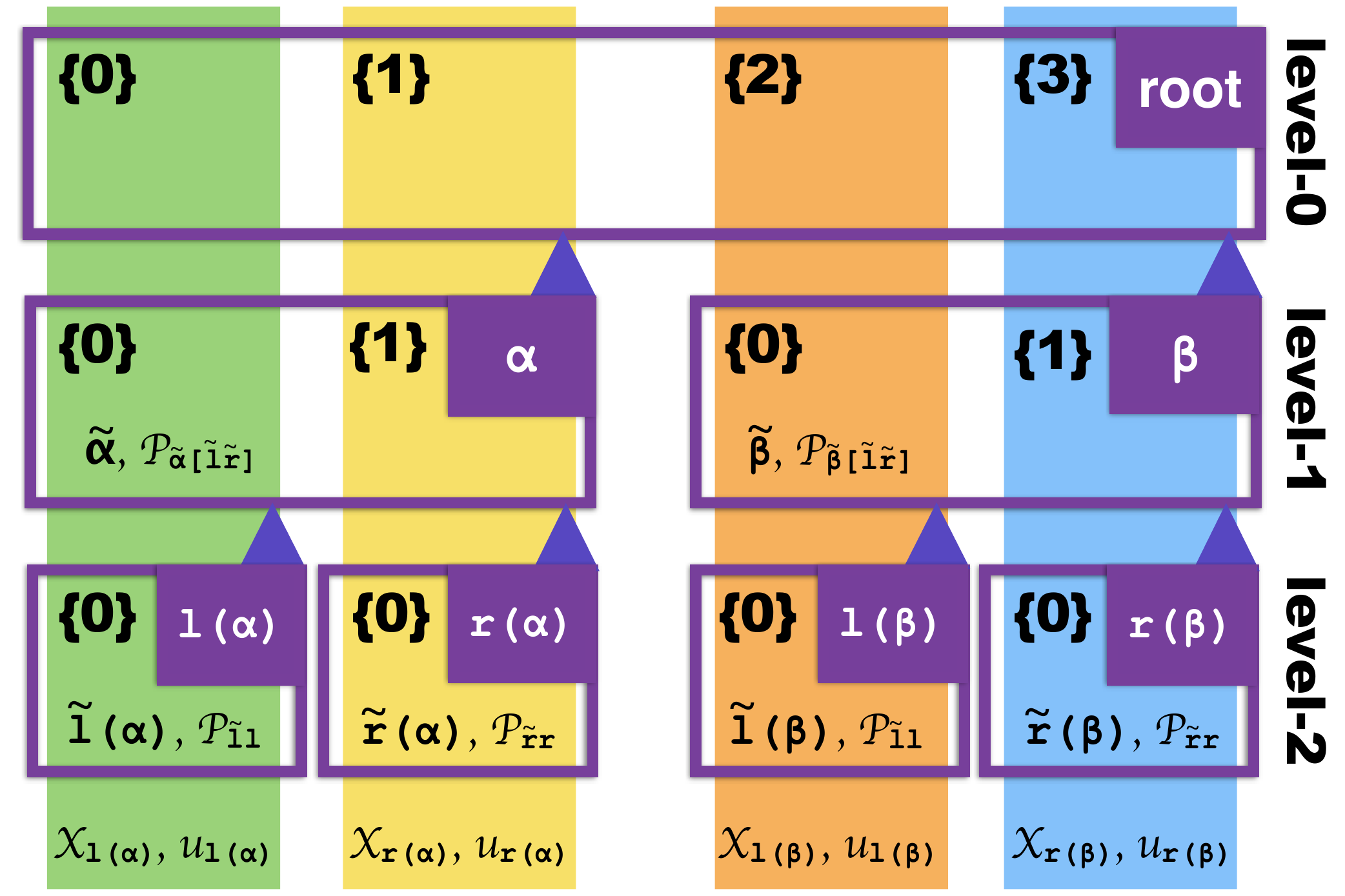} \caption{The \ASKIT{}
  distributed binary tree with $p=4$.  The purple boxes indicate tree
  nodes.  Depending on the level a tree node is logically distributed
  across several MPI processes, which are organized under an
  MPI-communicator. The ranks are indicated in curly brackets.  Each
  process owns the data of its colored rectangle region.  For example,
  the yellow process owns $\rc(\alpha)$ and $P_{\tilde{\rc}\rc}$ at
  level 2, but it doesn't own $\tilde{\alpha}$ at level 1. We depict
  the main data used in the factorization, the matrices $P$ and the
  skeleton points at each node. Notice that the nodes at level 2 are
  not the leaf nodes. But the computations in the subtrees at higher
  levels are entirely local to the MPI processes that own them. In
  this example we have four such processes, indicated by different
  colors.} \label{fig:tree}
\end{figure}

In~\cite{march-xiao-yu-biros15} we discussed the parallelization
of \ASKIT{}. Here we build on that work and derive our parallel
factorization algorithm. We assume that we have $p$ ranks in the
global MPI communicator. The data-structure used in \ASKIT{} is
summarized in \figref{fig:tree}. The distributed-memory parallelism in
the algorithm is invoked between tree levels-0 and -$\log p$.  We
briefly review the structure of the tree before discussing the
algorithms and the complexity.

\textbf{Distributed Tree.}
Without loss of generality, we assume $N$ and $p$ are powers of two
and that each MPI process owns $n=N/p$ points.
\figref{fig:tree} is a distributed tree with $p=4$.
We use purple boxes to indicate tree nodes.  Each treenode node is
assigned
%
%\footnote{Here ``assigned'' means that $\bigO(1)$ data pertaining to
%a node at level $l$ are replicated and $\bigO(N/2^l)$ data are
%distributed across the processes in the node's communicator.}
%
to a number of MPI processes.  All processes sharing a node are
grouped to an MPI communicator. We use \{i\} to indicate the
$i_{\mathrm{th}}$ MPI rank in the local communicator.  For example, at
level 0 the root of the tree is distributed among all 4 processes
(from green to blue in~\figref{fig:tree}).  The green process has
local MPI rank \{0\}, and the blue process has rank \{3\}.  If a node
$\alpha$ contains $q$ processes, then the processes with rank
\{$i<\frac{q}{2}$\} are assigned to the left child of $\alpha$, and
processes with rank \{$i\geq\frac{q}{2}$\} are assigned to the right
child of $\alpha$. At level $\log{p}$, each tree node is assigned to a
single MPI process, which is also assigned the entire subtree of that
node from level $\log{p}+1$ to the leaf level (not displayed
in \figref{fig:tree}).

\begin{figure}[!t]
  \centering \includegraphics[scale=.4]{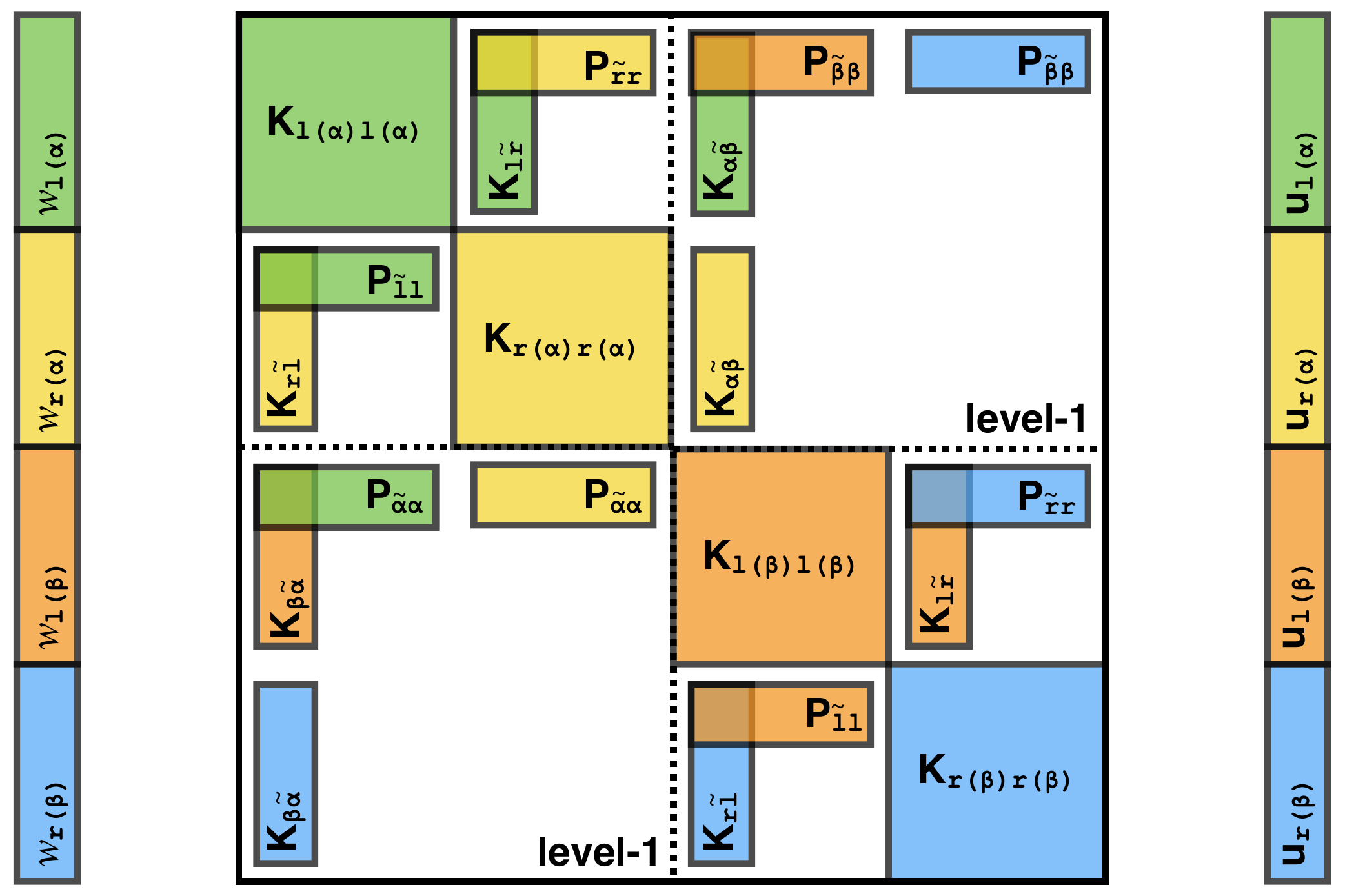} \caption{Distributed
  $K_{\alpha\sk{\beta}}$ and $P_{\alpha\sk{\alpha}}$ factors of the
  $\sk{K}$ matrix that corresponds to the partitioning
  of~\figref{fig:tree}. Each MPI processes (in the global
  communicator) stores only the matrices that have the same color with
  it. A matrix that has multiple color is distributed across several
  processes. We also show the
  partitioning of input and output vectors.}  \label{fig:matrix}
\end{figure}

Some processes may own the skeletons and $P$ matrix of the distributed
node. We give an example in \figref{fig:tree}; only \{0\} in each node
$\alpha$ stores the skeletons $\sk{\alpha}$ and
$P_{\sk{\alpha}[\lc\rc]}$. In the factorization we need to build
matrices $U_\alpha,W_\alpha,Z_\alpha$, which in turn involve
$P_{\sk{\lc}\lc}$, $P_{\sk{\rc}\rc}$, $K_{\lc\sk{\rc}}$ and
$K_{\rc\sk{\lc}}$. These matrices are unnecessary in \ASKIT{} but
in \IASKIT{} we need to compute and store them in parallel. These
computations involve collective broadcasts and reductions in the
communicators of $\lc,\rc$ and point-to-point exchanges between
the \{0\} and \{$q/2$\} ranks of the processes that were assigned to
$\alpha$.  Finally matrix $Z$ is stored at \{0\}. \figref{fig:matrix}
shows how these matrices are distributed for the matrix generated by
the tree in~\figref{fig:tree}.
%For simplicity,
%we use $U^{(l)}$ and $V^{(l)}$ to indicate the $UV$ factor in some
%tree nodes at level $l$.  

%Despite that the tree nodes are replicated, but the data the
%node contains is distributed. Other than the data point $\XX$
%and potential vector $\bu$, processes may possess the skeletons
%information.
%For example, the orange process
%possesses $\rc(\alpha)$ and $P_{\tilde{\rc}\rc}$ at level 2,
%but it doesn't possess $\tilde{\alpha}$ in $alpha $at level 1.
%In brief, only \{0\} in each node possess the skeletons and $P$.

%\textbf{Distributed matrices.}
%Before we go into parallel algorithms, we first explain how
%$K_{\alpha\sk{\beta}}$ (overwritten by $W_{\alpha\sk{\beta}}$)
%and $P_{\alpha\sk{\alpha}}$ are distributed. 
%Recall that for each node $\alpha$, we approximate $K$ in
%\eqref{e:blocksmw} and \eqref{e:innerid}.
%Since $\alpha$ is distributed, so does 
%$K_{\alpha\sk{\beta}}$ (overwritten by $W_{\alpha\sk{\beta}}$)
%$V_{\alpha}$. 
%\figref{fig:matrix} shows how these factors are distributed.
%For simplicity, we use $U^{(l)}$ and $V^{(l)}$ to indicate the $UV$ 
%factor in some tree nodes at level $l$. 
%A process owns only portions with the same color.

\begin{algorithm}[!htp]
\begin{algorithmic}[1]
  \IF {$\alpha$ is at level $\log{p}$}
    \STATE \texttt{Factorize($\alpha$)}.
  \ELSE
    \STATE \{$<\frac{q}{2}$\} $\texttt{DistributedFactorize}(\lc, \frac{q}{2})$.
    \STATE \{$\geq\frac{q}{2}$\} $\texttt{DistributedFactorize}(\rc, \frac{q}{2})$.
    \STATE \{0\} \texttt{Bcast} $P_{\tilde{\alpha}[\tilde{l}\tilde{r}]}$ inside the distributed node.
    \STATE \{0\}$\rightarrow$\{$\frac{q}{2}$\} \texttt{SendRecv} $\sk{\lc}$.
    \STATE \{$\frac{q}{2}$\} \texttt{Bcast} $\sk{\lc}$ using the $\rc$ communicator.
    \STATE \{$\frac{q}{2}$\}$\rightarrow$\{0\} \texttt{SendRecv} $\sk{\rc}$.
    \STATE \{0\} \texttt{Bcast} $\sk{\rc}$ using the $\lc$ communicator.
    \STATE \{$<\frac{q}{2}$\} Compute $K_{\lc\sk{\rc}}$, $P_{\sk{\lc}\lc}$ portions.
    \STATE \{$\geq\frac{q}{2}$\} Compute $K_{\rc\sk{\lc}}$, $P_{\sk{\rc}\rc}$ portions.
    \STATE \{$<\frac{q}{2}$\}$W_{\lc\sk{\rc}}=\texttt{DistributedSolve}(\lc, K_{\lc\sk{\rc}}, \frac{q}{2})$.
    \STATE \{$\geq\frac{q}{2}$\}$W_{\rc\sk{\lc}}=\texttt{DistributedSolve}(\rc, K_{\rc\sk{\lc}}, \frac{q}{2})$.
    \STATE \{0\} \texttt{Reduce} $P_{\sk{\lc}\lc}W_{\lc\sk{\rc}}$.
    \STATE \{$\frac{q}{2}$\} \texttt{Reduce} $P_{\sk{\rc}\rc}W_{\rc\sk{\lc}}$.
    \STATE \{$\frac{q}{2}$\}$\rightarrow$\{0\} \texttt{SendRecv} $P_{\sk{\rc}\rc}W_{\rc\sk{\lc}}$.
    \STATE \{0\} computes and LU factorizes $Z$.
  \ENDIF
\end{algorithmic}
\caption{{} \texttt{DistributedFactorize}($\alpha$,$q$)}
\label{a:distributedfactor}
\end{algorithm}

\textbf{Factorization.}
In \algref{a:distributedfactor}, we describe how \eqref{e:blocksmw}
and \eqref{e:innerid} are computed.  Each process uses its local MPI
rank to make different decisions. E.g., if a statement is labeled with
$\{<\frac{q}{2}\}$, then only processes with rank $<\frac{q}{2}$
execute it.  Otherwise, all statements are executed by all processes.
The whole communication pattern is simple. If a process owns the
skeletons and $P$ of $\alpha$, then it sends $P$ to all processes in
the same node and sends its skeletons to all processes in its sibling
tree node.  In the pseudo-code, this is a little bit complicate, since
we can't simply broadcast between communicators. Thus, some
point-to-point communication is involved.
%When two sibling nodes
%communicate they use their root processes (in their own communicator)
%and the communicator of their parent.

We will explain the algorithm using the yellow process
(\textbf{Yellow}) in \figref{fig:tree} and \figref{fig:matrix}. 
At level 2, \algref{a:factor} is invoked, and the local tree of 
$\rc(\alpha)$ is factorized.
At level 1, \textbf{Yellow} returns from line 5 as rank \{1\}.
At line 6, it receives $P_{\alpha[\sk{\lc}\sk{\rc}]}$, which will 
later be used to compute $P_{\alpha\sk{\alpha}}$ in the parent level.
At line 7, it receives $\sk{\lc}$ from \textbf{Green}, and at line 9
it sends $\sk{\rc}$ it owns to \textbf{Green}.
$\sk{\lc}$ is used to evaluate $K_{\rc\sk{\lc}}$ at line 12.
Evaluating $P_{\sk{\rc}\rc}$ doesn't require communication at this level, 
since all factors in the subtree are owned by \textbf{Yellow}.
At line 14, $K_{\rc\sk{\lc}}$ is overwritten by $W_{\rc\sk{\lc}}$.
Finally, \textbf{Yellow} sends its $P_{\sk{\rc}\rc}W_{\rc\sk{\lc}}$ to 
\textbf{Green}, and $Z$ is computed and factorized on \textbf{Green}.
At level-0, we only illustrate how $P_{\sk{\alpha}}[\sk{\lc}\sk{\rc}]$
that 
\textbf{Yellow} received at level 1 is used to compute the $P_{\sk{\alpha}\alpha}^{yellow}$
at line 11. By the definition of \eqref{e:innerid},
the right portion of $P_{\sk{\alpha}\alpha}=[P_{\sk{\alpha}\alpha}^{green},
  P_{\sk{\alpha}\alpha}^{yellow}]$ only requires
  $P_{\sk{\alpha}[\sk{\lc}\sk{\rc}]}$ and $P_{\sk{\rc}\rc}$.
  The first part
  $P_{\sk{\alpha}[\sk{\lc}\sk{\rc}]}$ is received at level 1, and the second
  part $P_{\sk{\rc}\rc}$ is computed at level 1.

\begin{algorithm}[!htp]
\begin{algorithmic}[1]
  \IF {$\alpha$ is at level $\log{p}$}
    \STATE $u=\texttt{Solve}(\alpha, \bu)$. 
  \ELSE
    %\IF {local MPI rank $\{i\}<\frac{q}{2}$}
      \STATE $\{<\frac{q}{2}\}$ $u=\texttt{DistributedSolve}(\lc, \bu, \frac{q}{2})$.
    %\ELSE
      \STATE $\{\geq\frac{q}{2}\}$ $u=\texttt{DistributedSolve}(\rc, \bu, \frac{q}{2})$.
    %\ENDIF
    \STATE Compute $P_{\sk{\alpha}\alpha}\bu$ portion.
    \STATE \{0\} \texttt{Reduce} $P_{\sk{\alpha}\alpha}\bu$ portion and compute $ZP_{\sk{\alpha}\alpha}\bu$.
    \STATE \{0\} \texttt{Bcast} $ZP_{\sk{\alpha}\alpha}\bu$.
    \STATE Compute $\bw=\bu-W_{\alpha\sk{\alpha}}ZP_{\sk{\alpha}\alpha}\bu$.
  \ENDIF
\end{algorithmic}
\caption{{} $\bw=\texttt{DistributedSolve}(\alpha, \bu, q)$}
\label{a:distributedsolve}
\end{algorithm}

{\bf Solver.}  \algref{a:distributedsolve} uses a similar bottom-up tree
traversal. \algref{a:solve} is invoked when the traversal reaches
level-$\log{p}$.  At line 7, \{0\} reduces all portions of
$P_{\sk{\alpha}\alpha}\bu$ and computes $ZP_{\sk{\alpha}\alpha}\bu$
where $Z$ was previously factorized.  At line 8,
$ZP_{\sk{\alpha}\alpha}\bu$ is broadcast by \{0\}.  Finally, each
process computes its own
$\bw=\bu-W_{\alpha\sk{\alpha}}ZP_{\sk{\alpha}\alpha}\bu$.

{\bf Complexity.}  We use a network model to derive the parallel
complexity.  We assume that the network topology provides
$\MA{O}(\ns \log{p})$ complexity for \texttt{Bcast}
and \texttt{Reduce} for a message of size $\ns$. For notational
clarity we omit the latency cost.

Let $T^{s}(q)$ be the parallel complexity of \algref{a:distributedsolve}. 
$T^{s}(1)=\MA{O}(n\log{n})$ is the base case complexity of \algref{a:solve}.
We can derive the complexity recursively as
%\begin{equation}
$T^{s}(q) = T^{s}\left(\frac{q}{2}\right) + \MA{O}(ns+s^2) + \MA{O}(s\log{q}).$
%\end{equation}
Here $\MA{O}{(ns+s^2)}$ is the message transfer computation cost of \eqref{e:blocksmw},
and $\MA{O}(s\log{q})$ is the communication cost.  The overall \textbf{\emph{solve
complexity}} with grain size $n=N/p$ is given by
\begin{equation}\label{e:t-solv}
  T^s(n,p) = \MA{O}(n\log{N})_{comp} + \MA{O}(\ns\log^{2}{p})_{comm},
\end{equation}
including both computation and communication.  The number of messages
 required is $\MA{O}(\log^2{p})$ if one wants to include latency
 costs. We observe that the communication cost doesn't scale with the
 problem size $N$, but only with $s$ and $\log^2 p$.
%$\MA{O}(rn\log{n})$ is the complexity of solving the local subtree. 
%$\MA{O}(rn\log{p})$ is the distributed computation cost, $\MA{O}(rn)$
%per level. 
%$\MA{O}(r\log^{2}{p})$ is the communication complexity.
%Each level contains $\MA{O}(r\log{p})$ collective communication,
%and there are $\log{p}$ levels in total.

The complexity  $T^{f}(q)$  of \algref{a:distributedfactor} is
similarly given by
%\begin{equation}
$  T^f(q) = T^f\left(\frac{q}{2}\right) + s T^s\left(\frac{q}{2}\right) + \MA{O}(ns) + \MA{O}(sd+s^2).$
%\end{equation}
Here $sT^s(\frac{q}{2})$ is the time
in \algref{a:distributedsolve}, $\MA{O}(ns)$ the time
in \eqref{e:innerid}, and $\MA{O}(sd+s^2)$ is the communication cost
spent on broadcasting skeletons and $P$.  Overall, the \textbf{\emph{factorization
complexity}} with grain size $n=N/p$ is 
\begin{equation}\label{e:t-fact}
  T^f(n,p) = \MA{O}(n\log^{2}{N})_{comp} + 
             \MA{O}(s^2\log^{3}{p}+d\log{p})_{comm},
\end{equation}
 including both computation and communication.  Observe that the
sequential complexity is $p$ times the parallel computation
complexity, which indicates that the algorithm should be able to
achieve almost linear scalability.  The communication cost is minor
since it does not depend on $N$.

{\bf Discussion.}
%Despite that the communication cost is not a bottleneck, 
\IASKIT{} shares the common problems of all treecodes: (1) the
parallelism diminishes exponentially when traversing to the root, and
(2) the possible global synchronization during the level by level
traversal.  Solving \eqref{e:smw} in a distributed node reduces the
$PW$ product to one process in the communicator. This takes the
advantage of the collective \texttt{Bcast} and \texttt{Reduce}, but
the parallelism diminishes.  This efficiency degradation is not
significant when $p$ and $s$ are small, but we will observe this
effect in the scaling experiments.  The load-balancing problem happens
when the adaptive skeletonization is used to dynamically decide the
rank of the off-diagonal blocks.  The processes with the largest rank
becomes the bottleneck. This load-balancing issues requires a dynamic
scheduling and beyond the scope of this paper.

\section{Experimental Setup} \label{s:setup} We performed numerical experiments to examine the accuracy and the
scalability of \IASKIT{}. Notice that the only alternative method that
we could compare with \IASKIT{} for the datasets of the size and
dimensionality is a global low-rank approximation of $K$. But
in~\cite{march-xiao-biros-e15}, we have already shown that for several
bandwidths of practical interest the low-rank approximation is not a
good scheme for approximating $K$.

We also applied \IASKIT{} to  binary classification tasks on
real-world datasets using kernel regression. Next we summarize the hardware,
datasets, and parameters used in our experiments. In the majority of
our tests we use the Gaussian kernel function~\eqref{e:gaussian}.

\textbf{Implementation and hardware.}
All experiments were conducted on the Maverick and Stampede clusters at
the Texas Advanced Computing Center.  Maverick nodes have two 10-core
Intel Xeon E5-2680 v2 (3.1GHz) processors and 256GB RAM.  Stampede
nodes have two 8-core E5-2680 (2.8GHz) processors and 32GB RAM.
\IASKIT{} and \ASKIT{} are written in \texttt{C++}, employing OpenMP
for shared memory parallelism, Intel MKL for high performance linear
algebra operations and Intel MPI for distributed memory parallelism.
A custom direct ``MatVec''\footnote{See the General Stride Kernel
Summation project: \url{https://github.com/ChenhanYu/ks}.} in \ASKIT{}
is optimized with \texttt{SSE2} and \texttt{AVX}.  The classification
code is written in \texttt{C++}, employing a Krylov subspace method
(\GMRES{}) from the PETSc library~\cite{petsc-web-page}.

\textbf{Iterative methods.} 
Recall that \IASKIT{} can only compute $(\lambda I + \sk{K})^{-1}\bu$
without the sparse matrix $S$.  While $S \neq 0$, \ASKIT{} can better
approximate $K$ using nearest-neighbor pruning ($\epsilon_M$ will be
smaller), but we can no longer compute its exact inverse directly.
Still \IASKIT{} can be a good preconditioner.  In this case,
converging to the desired accuracy may take several GMRES iterations,
but the classification accuracy may also be improved.

\begin{table}
\centering
{ \begin{tabular}{rrrrrrrrr} \toprule Dataset && $N_{\mathrm{train}}$
  && $N_{\mathrm{test}}$ && $d$ &&
  $h$\\ \cline{1-1} \cline{3-3} \cline{5-5} \cline{7-7} \cline{9-9} \cellcolor{gray!50} \textbf{COVTYPE}
  && 500K && 10K && 54 && 0.07 \\ \cellcolor{gray!50} \textbf{SUSY} &&
  4.5M && 10K && 8 && 0.07 \\ \cellcolor{gray!50} \textbf{MNIST} &&
  1.6M && 10K && 784 && 0.10 \\ \cellcolor{gray!50} \textbf{NORMAL} &&
  1M-32M && - && 64 && 0.31-0.17 \\ \toprule \end{tabular}
  } \caption{Datasets used in the experiments. Here
  $N_{\mathrm{train}}$ denotes the size of the training set and
  corresponds to the $N$ in \IASKIT{}; $N_{\mathrm{test}}$ is the
  number of test points used to check the quality of the
  classification. $d$ is the dimensionality of points in the dataset
  and $h$ is the bandwidth of the Gaussian kernel used in our
  experiments. We also performed tests with a polynomial kernel.}
\label{tab:dataset}
\end{table}

\textbf{Datasets:}
In \tabref{tab:dataset}, we give details of each dataset.  In strong
scaling and the classification tasks on Maverick, we use three
real-world datasets: \textbf{SUSY} (high-energy physics)
and \textbf{COVTYPE} (forest cartographic variables) from 
%UCI machine learning repository 
\cite{Lichman:2013}, \textbf{MNIST} (handwritten digit recognition)
from \cite{chang2011libsvm}.  In the weak scaling tasks, we generate a
random 64D point cloud which are drawn from a 6D \textbf{NORMAL}
distribution and embedded in 64D with additional noise. This resembles
a dataset with a high ambient but relatively small intrinsic
dimension, something typical in data analysis.  We also report strong
scaling results using a 16M \textbf{NORMAL} synthetic dataset.

\textbf{Bandwidth selection and accuracy metrics.}
For the majority of experiments, we used the Gaussian kernel. The
values are taken from~\cite{march-xiao-biros-e15}, where they were
selected using cross-validation with a kernel Bayes classifier. The
regularization $\lambda$ is also selected by cross-validation. 

We test both versions of \ASKIT{} (defined by the point-node
separation criterion discussed in~\secref{s:askit}), which correspond
to $S=0$ and $S \neq 0$ in~\eqref{e:askit}. We use three metrics to
measure the quality of the approximations for the {\sc MatVec} and {\sc
Invert} approximations of \ASKIT{} and \IASKIT{}:%
\footnote{Reminder on the notation: Here $K$ is the exact kernel
matrix and $\sk{K}$ is its \ASKIT{} approximation with $S=0$. 
By $(\lambda I + \sk{K})^{-1}$, we refer to the inverse based on the 
\ASKIT{} approximation.}
\begin{equation}\label{e:error}
\begin{split}
\epsilon_M &= \|K\bw-(\sk{K}+S)\bw\|_2/\|K \bw\|_2, \\
\epsilon_I &= \|\bw-(\lambda I +\sk{K})^{-1}(\lambda I + \sk{K})\bw\|_2/\|\bw\|_2,\\
\rho  &= \| \bu-(\lambda I + \sk{K}+S)\bw)\|_2/ \|\bu\|_2.\\
\end{split}
\end{equation}
The first error measures the \ASKIT{} {\sc MatVec} error.  
If the off-diagonal blocks of $K$ indeed admit low-rank,
then $\epsilon_M$ should be small.
Since evaluating $K\bw$ exactly ($\MA{O}(N^2)$ work) is impossible for large $N$, 
we estimate $\epsilon_M$ by sampling 1K points.  
We measure the condition of \IASKIT{}'s {\sc Invert} for the case $S=0$ using $\epsilon_I$.  
If $(\lambda I + \sk{K})$ is well-conditioned, then $\epsilon_I$
should be small.
We test $(\lambda I + \sk{K})^{-1}$ as a
preconditioner for \GMRES{} for solving $(\lambda I + \sk{K}+S)\bw
= \bu$ and we report the normalized residual $\rho$ upon termination  
of \GMRES{}. 
Typically the desired $\rho$ for the classification 
tasks is 1e-3. However, if the \GMRES{} terminates while reaching the maximum
iteration steps, then $\rho$ will be larger. 
Overall, we measure the quality of a preconditioner according to the
required iteration steps \#iter and the achieved accuracy $\rho$. 

{\bf Kernel regression for binary classification.}  Given a set of
points with binary labels%
\footnote{The sets {\bf MNIST} and {\bf COVTYPE} have several classes. For those we
perform a one-vs-all binary classification.  Take {\bf MNIST} dataset
as an example, the two classes in the binary classification task are
whether a point is the digit $3$ or not.}
, we split these points into the training and
the testing sets. We use the training set to compute the weights for
the regression. That is, given the correct classification labels
$\bu$, where $\bu(i) \in \{-1,1\}$ for each training point $i$, we
solve $(\lambda I + \sk{K}+S)\bw = \bu$ for the regression weights
$\bw$. Once we compute these weights, the label for a test point is
given by $\bx \notin \XX$ is $\text{sign}(\Ker(\bx,\XX)\bw)$.  By
comparing the true labels of the test points (which are given) to the
computed ones using the kernel regression, we can estimate the
classification accuracy as a percentage of correctly classified
points.  The classification errors are reported in the next section.
%% The
 %% dataset is excluded from the training dataset, and the accuracy
%% averages over 10K instances.

\section{Empirical Results} \label{s:results} We present three sets of results, for weak and strong scaling and for
 kernel regression. We conducted 29 experiments labeled from \#1
 to \#29 in the tables.  In \tabref{tab:strong}, we report numerical
 accuracy and strong scaling results.  In \tabref{tab:weak}, we report
 weak scaling for both Gaussian and polynomial kernels.  Finally, we
 apply \IASKIT{} as a preconditioner in the binary classification
 problem.  We show that \IASKIT{}'s $(\lambda I + \sk{K})^{-1}$ speeds up the
 convergence rate up to 5$\times$ and thus it is a good preconditioner
 for $(\lambda I + \sk{K}+S)$. We expect this to be the case when $\|S\|$ is
 sufficiently smaller than $\sk{K}$. In the results the parameter
 $\tau$ indicates the accuracy tolerance used in the skeletonization
 phase of the \ASKIT{} algorithm.
 $s_{max}$ is the maximum skeleton size. $\kappa$ is the number
 of the \emph{nearest-neighbors}.

\begin{table}[!t]
\centering
{
  \begin{tabular}{r|r|r|rrrr|r} 
  \hline
  \# & core & $T^{M}$ & $T^{f}_{c}$ & $T^{f}_{f}$ & $T^{s}$ & $T$ & $\eta$ \\
  \hline
  \rowcolor{gray!50}
  \multicolumn{3}{l}{\textbf{COVTYPE}} & 
  \multicolumn{5}{l}{$h=0.07$, $\lambda=0.3$, $\kappa=2048$, $\tau=10^{-7}$,} \\
  \rowcolor{gray!50}
  \multicolumn{3}{l}{} & \multicolumn{5}{l}{$\epsilon_M=1e-1$, $\epsilon_I=4e-12$.} \\
  \hline
  \#1 & 160 & 1.1 &  $0.2$ & $80+64$ & $1.0$ &  145 & 1.00 \\ % 1.6 TFLOPS
  \#2 & 320 & 0.6 &  $0.4$ & $30+46$ & $0.8$ &   77 & 0.94 \\ % 2.9 TFLOPS
  \#3 & 640 & 0.3 &  $0.7$ & $10+31$ & $0.5$ &   42 & 0.86 \\ % 5.4 TFLOPS
  \hline
  \rowcolor{gray!50}
  \multicolumn{3}{l}{\textbf{SUSY}} & 
  \multicolumn{5}{l}{$h=0.07$, $\lambda=10$, $\kappa=2048$, $\tau=10^{-3}$,} \\
  \rowcolor{gray!50}
  \multicolumn{3}{l}{} & \multicolumn{5}{l}{$\epsilon_M=2e-1$, $\epsilon_I=3e-14$} \\
  \hline
  \#4 & 320 & 2.6 &  $1.2$ & $636+666$ & $3.8$ & 1307 & 1.00  \\
  \#5 & 640 & 1.3 &  $1.5$ & $269+409$ & $2.3$ &  681 & 0.96  \\ % 7.02 TFLOPS 
  \hline
  \rowcolor{gray!50}
  \multicolumn{3}{l}{\textbf{MNIST}} &
  \multicolumn{5}{l}{$h=0.3$, $\lambda=0$, $\kappa=256$, $\tau=10^{-3}$,} \\ 
  \rowcolor{gray!50}
  \multicolumn{3}{l}{} &\multicolumn{5}{l}{$\epsilon_M=6e-8$, $\epsilon_I=4e-13$} \\
  \hline
  \#6 & 160 &  2.5 &  $0.1$ & $12.2+3.0$ &  $0.26$ & 16 & 1.00 \\ % 0.37
  \#7 & 320 &  1.3 &  $0.1$ &  $6.6+2.1$ &  $0.15$ &  9 & 0.88 \\ % 0.66
  \#8 & 640 &  0.5 &  $0.1$ &  $3.6+1.6$ &  $0.08$ &  5 & 0.80 \\ % 1.19
  \hline
  \rowcolor{gray!50}
  \multicolumn{3}{l}{\textbf{NORMAL}} &
  \multicolumn{5}{l}{$h=0.19$, $\lambda=0.1$, $\kappa=128$,} \\ 
  \rowcolor{gray!50}
  \multicolumn{3}{l}{} &\multicolumn{5}{l}{$\epsilon_M=4e-1$, $\epsilon_I=4e-6$} \\
  \#9  & $1,024$ &   1.5 &  $0.2$ &  $37.6+57.3$ & $0.9$ & 96 & 1.00 \\
  \#10 & $2,048$ &   0.8 &  $0.3$ &  $17.4+33.7$ & $0.5$ & 52 & 0.92 \\
  \#11 & $4,096$ &   0.4 &  $0.2$ &  $8.6+21.0$ &  $0.3$ & 30 & 0.80 \\
  \hline
  \end{tabular}
}
\caption{Strong scaling results. Experiments \#1--\#8 are done
  on Maverick using $m=2,048$ and $s_{max}=2,048$.
  Experiments \#9--\#11 are done on Stampede with $N=16M$, $m=512$ and
  $s=256$.  $T^{M}$ is the {\sc MatVec} time, $T^f_f$ is computation
  time during the factorization, $T^f_c$ is communication time during
  factorization, and $T^s$ is solve time, all in seconds; $\eta$
  denotes the parallel efficiency. }
\label{tab:strong}
\end{table}

\textbf{Accuracy.} In~\tabref{tab:strong} and \tabref{tab:weak}, 
we report \ASKIT{}'s $\epsilon_M$ error for
  reference but our focus is on $\epsilon_I$ and $\rho$, which are
 related to \IASKIT{}.  Note that $\epsilon_I$ is not always small due
 to the amplification of rounding errors by the conditioning of
 $(\lambda I + \sk{K})$.  Both $h$ and $\lambda$ have strong impact on the condition
 number.  For a Gaussian kernel matrix, a smaller $h$ usually results
 in better conditioned systems since it approaches the identity
 matrix. Runs \#6--\#8 are examples where the systems are
well-conditioned even with $\lambda=0$.  With a larger $h$ and lower
$d$, $\tilde{K}$ becomes nearly rank-deficient.  In \#12--\#19 we see
that $\lambda=0.1$ is not sufficiently large, and that's why we observe
large $\epsilon_I$.
%The expected $\rho$ upon GMRES termination is 1e-3. \#20 and \#22
%fail to converge to this accuracy because of reaching the maximum 
%iteration steps.
%Typically the desired $\rho$ for the classification 
%tasks is 1e-3. However, if the \GMRES{} terminates while reaching the maximum
%iteration steps, then $\rho$ will be larger. 
%Overall, we measure the quality of a preconditioner according to the
%required iteration steps \#iter and the achieved accuracy $\rho$. 

\begin{table}[!t]
\centering
{\footnotesize \begin{tabular}{|l|rrrr|rrrr|} \hline
  & \multicolumn{4}{c|}{Gaussian} & \multicolumn{4}{c|}{Polynomial} \\ 
  \hline \# & \#12 & \#13 & \#14 & \#15 & \#16 & \#17 & \#18 & \#19 \\ 
  \hline 
  $N$ & 1M & 4M & 16M & 32M & 1M & 4M & 16M & 32M \\ 
  {\bf cores} & 20 & 80 & 320 & 640 & 20 & 80 & 320 & 640 \\ 
  $h$ & 0.31 & 0.24 & 0.19 & 0.17 & - & - & - & - \\
  \hline
  $\epsilon_M$ & 1e-1 & 2e-1 & 3e-1 & 3e-1
  & \multicolumn{4}{c|}{\cellcolor{gray!50}$\MA{O}(1e-12)$} \\
  $\epsilon_I$ & 1e-7 & 1e-7 & 6e-7 & 1e-6 & 6e-8 & 3e-7 & 1e-6 &
  2e-6 \\
  %$K\bw$         &  3.3 &  5.4 &  5.5 &  5.6 &  96.1 & 149.7 & 165.0 &   165 \\
  \hline
  $T^{f}_{c}$  &  0.0 &  0.1 &  0.1 &  0.2 &   0.0 &   0.0 &   0.1 &   0.1 \\
  $T^{f}_{f}$  &  295 &  396 &  509 & 574  &   300 &   401 &   509 &   574 \\
  %$T^{s}_m$      &  0.0 &  0.0 &  0.0 & 0.0  &   0.0 &   0.0 &   0.0 &   0.0 \\
  $T^{s}$        &  2.5 &  3.1 &  3.8 & 4.2  &   2.5 &   3.1 &   3.8 &   4.2 \\
  $T$            &  298 &  399 &  513 & 578  &   303 &   404 &   513 &   578 \\
  \hline
  $\eta_f$    &  0.3  &  1.1 &  4.2 &  8.2 &   0.3 &   1.1 &   4.3 &   8.2 \\
  $\eta$      & 1.00  & 0.92 & 0.88 & 0.85 &  1.00 &  0.92 &  0.89 &  0.85 \\
  \hline
  \end{tabular}
}
\caption{Weak Scaling experiment results on \textbf{NORMAL} synthetic datasets.
  All experiments use $m=512$, $s=512$, $\kappa=128$ and
  $\lambda=0.1$.  The runtime is measured in seconds.  $\eta_f$ denotes
  \texttt{TFLOPS} and $\eta$ denotes the relative efficiency.}
\label{tab:weak}
\end{table}

\textbf{Scaling.}
Now, we discuss the weak and strong scaling results
in \tabref{tab:strong} and \tabref{tab:weak}.  The overall \IASKIT{}
time $T=T^{f}_{c}+T^{f}_{f}+T^{s}$ is divided into
\textbf{\emph{factorization time}} $T^{c}_{f}+T^{f}_{f}$ and
the \textbf{\emph{solving time}} $T^{s}$.  The factorization time is
split into (1) communication $T^{f}_{c}$ and (2) computation
$T^{f}_{f}$. $T^{f}_{f}$ is further subdivided into local and
distributed time. The communication part of $T^{s}$ is negligible so
we do not split $T^{s}$.  Take \#1 as an example, the communication
time is $0.2$ seconds, the local factorization time is $80$ seconds,
and the distributed factorization time is $64$ seconds. The solving
time is $1$ second, and the overall time is $145$ seconds.

In \tabref{tab:strong} we report the strong-scaling efficiency
$\eta$. The base-line experiments are \#1, \#4, \#6 and \#9, which have 
$\eta=1.00$. Fixing the problem size, we measure the efficiency 
by doubling the core number each time.
We observe that the communication cost $T^{f}_{c}$ is
insignificant; roughly speaking, doubling the number of cores will
reduce the efficiency by 4\%--12\%.  The efficiency degradation is
due to diminishing parallelism during the computation
of \eqref{e:smw}, where only the root rank in each subcommunicator
works.  We also measure the relative-to-peak FLOPS efficiency.  The
peak for 640-core on Maverick is $15.8$\footnote{\scriptsize{Each
E5-2680 v2 processor can compute $3.1\times8\times10=248$ GFLOPS per
second.  The aggregate theoretical peak performance of 64 processors
is $64\times248=15,782$ GFLOPS $\approx 15.8$ TFLOPS.}}  TFLOPS, and the peak of 4096-core on
Stampede is $85.2$ TFLOPS.  The ratio FLOPS$/$peak from
low-to-high is \textbf{MNIST2M} 8\%,
\textbf{NORMAL} 20\%, \textbf{COVTYPE} 36\% and 
\textbf{SUSY} 45\%.

The scalability varies between datasets during the factorization,
because there are load-balancing issues in the level-by-level
traversal.  Adaptive skeletonization results in different values of $s$ for
nodes in the same tree level.  The node with the largest $s$ becomes
the critical path in each level. The load-balancing problem is
especially serious in \textbf{MNIST} where the workload can vary by
orders of magnitude.  On the other hand, \textbf{SUSY} has the most
balanced workload, which allows it to reach 45\% peak floating point
efficiency.  We also report the {\sc MatVec} performance in the
$T^{M}$ column. The {\sc MatVec} in \ASKIT{} is highly optimized and
does not require any communication. Typically, it can reach
near-linear scaling and more than 80\% of peak floating point
efficiency.

In \tabref{tab:weak}, we report weak-scaling results.  We perform
experiments using Gaussian and polynomial kernels and we observe
similar behavior. We report the FLOPS rate ($\eta_{f}$) 
and relative parallel efficiency ($\eta$).
The base-line experiments are \#12 and \#16.
The efficiency drops dramatically from 20-core runs to 80-core runs, since \#12
and \#16 do not involve MPI.  From 80 cores to 640 cores, doubling the
number of cores results only in 3\% loss, which shows that the weak
scalability is almost linear.  The weak-scalability experiments can
reach a higher efficiency, since there are no communication and
load-balancing issues.  As a result, \IASKIT{} can reach 52\%
theoretical peak floating point efficiency.

%{\footnotesize
%  \begin{tabular}{crrrrrrrrrr}
%  \toprule
%  \multicolumn{9}{l}{$H^{-1}u$, Polynomial} \\
%  \toprule
%       && 1M && 4M && 16M && 32M \\
%  \cline{1-1} \cline{3-3} \cline{5-5} \cline{7-7} \cline{9-9}
%  core           &&    20 &&    80 &&   160 &&  320 \\
%  $\epsilon_M$ && 1e-14 && 4e-14 && 2e-14 && 1e-13 \\
%  $\epsilon_I$ &&  6e-8 &&  3e-7 &&  1e-6 && 2e-6 \\
%  $T^{Kw}$       &&  96.1 && 149.7 && 165.0 &&  165 \\
%  \toprule
%  $T^{f}_m$      &&   0.0 &&   2.1 &&  11.7 &&   23 \\
%  $T^{f}_c$      &&   301 &&   394 &&   501 &&  554 \\
%  $T^{s}_m$      &&   0.0 &&   0.2 &&   0.5 &&  0.7 \\
%  $T^{s}_c$      &&   2.5 &&   2.9 &&   3.3 &&  3.5 \\
%  \toprule
%  $T$            &&   303 &&   399 &&   517 &&  581 \\
%  Efficiency     &&  1.00 &&  0.92 &&  0.84 && 0.81 \\
%  \toprule 
%  \end{tabular}
%}
%\caption{Results for Uniform dataset, $\lambda=0.1$, $\kappa=128$
%  $m=512$, $s=512$.}
%\label{tab:polynomial}
%\end{table}

\begin{table}[!t]
\centering
{\footnotesize
  \begin{tabular}{rrrrrrr}
  \toprule 
  \# & $M$ & $T_{\mathrm{train}}$ & \#iter & Train (\%) & Test (\%) & $\rho$ \\
  \midrule 
  %\rowcolor{gray!50}
  \multicolumn{2}{l}{\textbf{COVTYPE}} & 
  \multicolumn{5}{r}{$(\lambda I + \tilde{K})\bw = \bu$.} \\
  \midrule 
  \#20 &                   $I$ & 231 & 100 & 99\% & 96\% &  4e-3 \\
  \#21 & $(\lambda I + \tilde{K})$  &    62 & 0 & 99\% & 96\% & 2e-11 \\
  \midrule
  %\rowcolor{gray!50}
  \multicolumn{2}{l}{\textbf{COVTYPE}} & 
  \multicolumn{5}{r}{$(\lambda I + \tilde{K}+S)\bw = \bu$.} \\
  \midrule 
  \#22 &                   $I$ & 1725 & 100 & 99\% & 97\% & 3e-3 \\
  \#23 & $(\lambda I + \tilde{K})$ &   950 & 43 & 99\% & 97\% & 1e-3 \\
  \midrule
  %\rowcolor{gray!50}
  \multicolumn{2}{l}{\textbf{SUSY}} & 
  \multicolumn{5}{r}{$(\lambda I + \tilde{K})\bw = \bu$.} \\
  \midrule 
  \#24 &                     $I$ & 3056 & 49 & 80\% & 79\% & 1e-3 \\
  \#25 & $(\lambda I + \tilde{K})$ &   870 & 0 & 80\% & 79\% & 1e-14 \\
  \midrule
  %\rowcolor{gray!50}
  \multicolumn{2}{l}{\textbf{SUSY}} & 
  \multicolumn{5}{r}{$(\lambda I + \tilde{K}+S)\bw = \bu$.} \\
  \midrule 
  \#26 &                     $I$ & 11519 & 47  & 80\% & 79\% & 1e-3 \\
  \#27 & $(\lambda I + \tilde{K})$ &   3419 & 7 & 80\% & 79\% & 1e-3 \\  
  \midrule
  %\rowcolor{gray!50}
  \multicolumn{2}{l}{\textbf{MNIST}} &
  \multicolumn{5}{r}{$(\lambda I + \tilde{K})\bw = \bu$.} \\
  \midrule 
  \#28 &                     $I$ &  185 & 5 & 100\% & 100\% & 1e-9 \\
  \#29 & $(\lambda I + \tilde{K})$ &   36 & 0 & 100\% & 100\% & 2e-14 \\
  \toprule 
  %\rowcolor{gray!50}
  %\multicolumn{10}{l}{\textbf{MNIST}, $\tilde{K}+E$, $h=0.1$, $\kappa=256$, $\tau=1e-3$.} \\
  %\toprule 
  %\#31 &         $I$ &&    s() && 100\% && 100\% &&    - \\
  %\#32 & $\tilde{K}$ &&    s() && 100\% && 100\% && 2e-2 \\
  %\toprule 
  \end{tabular}
}
\caption{Classification results using 640 cores. Other parameters are identical to
the experiments in \tabref{tab:strong}. Here $M$ denotes the
preconditioner, which can be either $I$ or $\tilde{K}$.  The number
if \GMRES{} iterations is denoted by \#iter. When it is $0$, it means
immediate convergence because \IASKIT{} is highly accurate.  The GMRES
solver will terminate while reaching 100 iteration steps or
$\rho<1e-3$. }
\label{tab:classification}
\end{table}

\textbf{Classification.}
We conduct three classification experiments,
following~\cite{march-xiao-biros-e15}. The training time
$T_{\mathrm{train}}$ is presented in pairs to show both the runtime
and the iteration steps of the \GMRES{} solver.  Four experiments are
conducted for each dataset, combining different kernel matrix
approximations ($\lambda I + \tilde{K}$ or $\lambda I + \tilde{K}+S$) and different
preconditioners ($I$ indicates no preconditioning or $\lambda I + \tilde{K}$).  We
present both training and testing (10K samples) accuracy, and we
present the normalized residual $\rho$ upon termination
of \GMRES{}.

 Approximating the kernel matrix using $\tilde{K}$, we reach 96\%
accuracy in \textbf{COVTYPE}, 79\% accuracy in \textbf{SUSY} and 100\%
accuracy in \textbf{MNIST}.  Using the nearest-neighbor-based
separation criterion which creates a nonzero $S$, we can reach a
higher accuracy (97\%) in \textbf{COVTYPE}.  In \#21, \#25 and \#29,
we use $(\lambda I + \tilde{K})^{-1}\bu$ as initial guess, which allows \GMRES{} to
converge immediately to a high precision by only applying the
preconditioner.  The training time only
involves the factorization time and the solving time, which are both
deterministic.  Without preconditioners, \GMRES{} takes 3$\times$ to
5$\times$ longer to converge to the desired accuracy. During the
cross-validation, time spent on training may be much longer
depending on the combination $h$ and $\lambda$.  

In runs \#23
and \#27, we show that \IASKIT{} also works well as a preconditioner
when we include the nearest-neighbors $S$ in the \ASKIT{} {\sc MatVec}
approximation. 
Comparing \#22 to \#23, we can find that \#23 with \IASKIT{} only takes 43
iterations to converge to 1e-3. 
On the other hand, \#22 reaches the maximum iteration limit, only
converging to 3e-3.
Involving nearest-neighbors $S$ can better approximate $K$ in some
cases. For example, \#23 can reach a higher classification accuracy
than \#21.
Overall \IASKIT{} provides a 2$\times$ to 3$\times$
speedup. For nonlinear classification methods like kernel logistic
regression and support vector machines, the factorization can be
amortized across several nonlinear iterations. In those cases the
speedups will be even more dramatic.

\section{Conclusions} \label{s:conclusion} We have introduced a new parallel fast factorization algorithm for
kernel matrices defined on datasets in high dimensions.  We evaluated
our methods on both real-world and synthetic datasets with different
kernel functions and parameters. We conducted analysis and experiments
to study the complexity and the scalability of \IASKIT{}.  These
experiments include scaling up to 4,096 cores, calculations on very
large data sets and achieving up 50\% of peak FLOPS performance.  We
show that applying \IASKIT{} to classification problems can result in
fast convergence and competitive accuracy.  We do not assume symmetry,
thus we automatically support variable bandwidth kernels. To our
knowledge, \IASKIT{} is the only algorithm that can be used with
high-dimensional data for a large variety of kernels and the only
algorithm that supports shared and distributed memory parallelism.

%Hierarchical matrices and their inverses are well-studied and 
%known to be very robust in low dimensional space. 
%We study high dimensional kernel machines to
%understand the limitations, and we exploit the spatial data %distribution not
%only with its low-rank structure but also the sparsity.
%We develop a fast inverse that scales $\MA{O}(N\log^{2}{N})$ to %approximate \eqref{e:primary}, and we analyze its accuracy and %scalability.
%We show that \IASKIT{} is not only fast but also sufficient 
%accurate in high dimensional kernel classification tasks. 
%Although there is no known algorithm to inverse $\tilde{K}+E$
%fast, but we show that \IASKIT{} can be a good preconditioner 
%for $\tilde{K}+\lambda I+E$ or even other approximation.

Our future work will focus on further optimization of our
implementation. In particular, there are opportunities to resolve the
parallelism shrinking problem by dynamic scheduling.  Also, we plan to
explore other possible variants with lower factorization cost,
extensions to GPU architectures, and integration with other
higher-level algorithms like support vector machines.

%\newpage
\section*{Acknowledgment} 
This material is based upon work supported by AFOSR grants
FA9550-12-10484 and FA9550-11-10339; by NSF grant CCF-1337393; by the
U.S. Department of Energy, Office of Science, Office of Advanced
Scientific Computing Research, Applied Mathematics program under Award
Numbers DE-SC0010518 and DE-SC0009286; by NIH grant 10042242; by DARPA
grant W911NF-115-2-0121; and by the Technische Universit\"{a}t
M\"{u}nchen---Institute for Advanced Study, funded by the German
Excellence Initiative (and the European Union Seventh Framework
Programme under grant agreement 291763).  Any opinions, findings, and
conclusions or recommendations expressed herein are those of the
authors and do not necessarily reflect the views of the AFOSR, the
DOE, the NIH, or the NSF. Computing time on the Texas Advanced
Computing Centers Stampede system was provided by an allocation from
TACC and the NSF.

\newpage
\bibliographystyle{siam}
\bibliography{gb,refs}

\end{document}